\documentclass[11pt]{article}

\usepackage{amsmath}
\usepackage{amssymb}
\usepackage{amsthm}

\usepackage{graphicx}

\usepackage{ifthen}

\newcommand{\limn}{\lim_{n\rightarrow \infty }}

\newcommand{\Rb}{{\mathbb{R}}}
\newcommand{\Cb}{{\mathbb{C}}}

\newcommand{\xast}{x^*}
\newcommand{\xn}{ (x_{n}) } 
\newcommand{\setan}{{}_{\stackrel{\textstyle \longrightarrow}{n}}}

\newcommand{\neu}{NEU_m  }
\newcommand{\hyp}{HYP_j  }

\newtheorem{lemma}{Lemma}
\newtheorem{proposition}{Proposition}

\newtheorem{definition}{Definition}

\theoremstyle{definition}

\begin{document}

  \title{ Double newtonisation of fixed point sequences}
 \author{M\'{a}rio M. Gra\c{c}a 
\thanks{ 
Department of Mathematics,  
Ins\-ti\-tuto Superior T\'{e}cnico,
Universidade T\'{e}cnica de Lisboa, Av. Rovisco Pais, 1049--001, 
Lisboa, Portugal, 
  mgraca@math.ist.utl.pt\ .  This work was 
 supported  by IDMEC--IST through 
the  program POCTI of FCT (Portugal).} }  
\maketitle

  \begin{abstract} A neutral fixed point of a real iteration map $u$ 
  becomes a  super attracting fixed point using a suitable double newtonisation. 
  The map $u$ is so transformed into a map $w$ which is here called the  
  standard  accelerator of $u$. The map $w$
  provides a unifying process to deal with a large set of 
    fixed point sequences which are not convergent or converge slowly. Several examples illustrate the main results 
  obtained. 
\end{abstract}

 

\section{Introduction}

The computation of a fixed point $\xast$ of a real iteration map $u$ 
is a difficult task when  $| u'(\xast) |\geq 1$. The case $| u'(\xast) | 
=1$ corresponds to a {\em neutral or nonhyperbolic} fixed point whereas for
$| u'(\xast) |> 1$ the fixed point is said to be {\em repelling}  \cite{Hol}  . The aim 
of this work is to give a simple method to deal with these kind of 
difficulties.

Let us denote  by $NEU_m$, or simply by $NEU$, 
 the set of real iteration maps $u$,   $m\,\,(m\geq 
1)$ times continuously differentiable in a neighbourhood of a {\em 
neutral}  
  fixed point $\xast$ such that $u'(\xast)=1$.  In this work we construct simple accelerators for a quite large subset of 
$NEU$. Although, as it is   well known,  the set  of logarithmically convergent sequences, 
 $Log$  ,
cannot be acce\-le\-rated \cite{del,delgb},
 the same simple accelerators are also useful to apply on
a large set of logarithmic 
fixed point sequences, denoted by $Log_m$.

In order to show that a certain subset of $NEU$ is accele\-rable we follow the  
strategy: for a given 
 $u\in \, NEU$,    (i)   obtain an iteration function 
$v$ (depending on $u$)
for which $\xast$ is an attracting or repelling fixed point and (ii) 
   combine $v$ with $u$ in order to get a final ite\-ra\-tion 
function $w$ for which $\xast$ is   a super attractor. We will call 
such a $w$ an  {\em accelerator} of $u$. Notice that the basic assumption 
is the existence of an isolated fixed point $\xast$ for $u$ which is (in general) 
unknown.

Examples of accelerators constructed using the steps (i) and 
(ii) above  are  the 
  so--called {\em combined iteration functions} which have been introduced 
  and studied by the author  in \cite{graca1,graca3}.
  Under mild assumptions 
  on $u\in \, NEU$ it is  possible to obtain  
  a combined hyperbolic ite\-ration function for instance $v=C(x,u)$ and finally an 
  accelerator $w=C(x,v)$. We call
   $w$   {\em standard accelerator}.
 
In propositions \ref{prop1} and \ref{prop2}   simple 
accelerators are given for particular iteration maps in $NEU$. The main results of 
this paper are propositions \ref{prop5}--\ref{prop8}. All results 
concerning the standard accelerator are new. In lemma  \ref{prop7} and 
proposition \ref{prop8}  we characterize the {\em kernel} of the standard 
accelerator giving a necessary and sufficient condition for a certain map 
$u$ in $NEU$ to belong to the kernel of $w$. An iteration map $u$ belongs 
to the kernel of a map $h$ if $h(x)=\xast$ for any $x$. This means that the
 fixed point $\xast$ of $u$ can be computed using only 
one iteration of $h$.

In section \ref{sub1} we apply the results obtained for the standard map 
to show that a certain set of logarithmic fixed point sequences, denoted by $FIX$, is 
accelerable.

In section  \ref{sec3}   we present several examples illustrating the 
properties of the accelerators studied in previous sections. In examples 1 
and 2   simple accelerators are constructed for well known iteration maps 
namely  the function $u(x)=\sin(x)$ and the logistic map \cite[pp. 2]{stuart1}.
Examples 3 and 4 deal with   certain iteration maps   which  attracted the attention of many
 researchers (see for instance 
\cite{sab1, sab2, sedogbo}) with the aim of finding accelerators. We show 
that such maps are accelerable or even  
 belong to 
the kernel of our standard accelerator.
Finally, in example 5, the standard accelerator is used for the 
computation of a multiple zero of  a 
complex function.

 This work extends the author's results presented in \cite{graca3}.

\section{Definitions}\label{secdef}

In this work we will consider two classes of ite\-ra\-tion  
functions  defined as follows.

\begin{definition}\label{def1}

Assume there exists  an isolated fixed point $\xast$ of a given iteration function 
$u$, and that $u$ is defined and $m$--fold continuously 
differentiable in a neighbourhood ${\cal D}$ of $\xast$.  Define
  $$\label{eqa}
\neu=\left\{
 u\in  C^m (D):\, u(\xast)=\xast, 
\,u'(\xast)=1
 \right\}, 
$$
and, for   $j\geq 1 $,
$$\label{eqb}
 \hyp (u)=\left\{ 
v\in  C^j (D):\,
 v(\xast)=\xast, \,
\, |v'(\xast)|\neq 1 
 \right\}.
$$
\end{definition}  

The set $\neu$ is formed by iteration functions with a  neutral  
  fixed point $\xast$ while $\hyp (u)$ is the set of all 
iteration functions $v$ possessing the same fixed point than $u$ for which $\xast$ is {\em  hyperbolic } 
(attracting or repelling). When the degree of  smoothness, $m$ or $j$, is 
implicit or not relevant   these sets  will be 
denoted respectively by $NEU$ and $HYP$.

Note that the  definition of $NEU$ does not include 
iteration functions such that $u'(\xast)=-1$. However this is not 
restrictive since    for 
  $\tilde u(x) =   u(u(x))$ we get $\tilde u'(\xast)=1$, 
so $\tilde u\in \,NEU$.

For a given $u\in \,NEU$   we are interested first of all in obtaining $v$ in 
$HYP$ ($v$ 
depending on $u$) such that $v'(\xast)\neq 1$, that is $\xast$ is a {\em 
hyperbolic} fixed point for $v$. Whenever we can explicitly 
construct  
an iteration function $h$, having the same fixed point as $u$ such that $h'(\xast)=0$
this   will be called an {\em accelerator} of 
$u$. When we are able to  
  construct an ite\-ra\-tion map $h $ which accelerates any 
iteration function $u$ in a certain subset of $NEU$  we say 
that such a set is 
{\em accelerable} and that  $h$ is an {\em accelerator} of this set.

  Whenever
an iteration function  $h$ accelerates $u\in\, NEU$ and there exists an integer $k$ such that
$ h'(\xast)=h^{(2)}(\xast)=\ldots = h^{(k)}(\xast) =0$ and 
$h^{(k+1)}(\xast)\neq 0$ 
we say that $h$ is a  {\em $k$--accelerator} of $u$. Also,   when the   derivatives of any order of $h$
 are zero at $\xast$, the 
iteration function $h$ will be called 
an {\em $\infty$--accelerator}.

\section{Simple accelerators}\label{sec2}

For a conveniently chosen  $x_0\in {\cal D}$, a  function $u$ in $NEU$ 
gives
rise to a sequence $x_{n+1}=u(x_n)$. If $\xn \setan \xast$, this sequence is 
known to be of   {\em 
logarithmic} convergence   since it satisfies $\limn 
(x_{n+1}-\xast)/(x_n-\xast)=1$ (see \cite[pp. 2]{Brezinski1}).

The next proposition shows that 
logarithmic fixed point sequences having zero as fixed point are 
ea\-sily accelerable.  
 \begin{proposition}\label{prop1}
  For any $g\in NEU_m  $ such that $\xast=0$ 
   the following iteration function $h_j$,
\begin{equation}\label{eq3} 
\begin{array}{l} 
 h_1 (x)=\int_0^{x} g(t)\, dt \\
{}\\
h_j(x)=\int_0^x h_{j-1} (t)\,dt,\,\,j=1,\ldots ,m
\end{array}
\end{equation}
 is at least a  {\rm(}$j-1${\rm)}--accelerator of $g$.
\end{proposition}
\begin{proof}
As $g\in NEU_m $ then $g(0)=0$ and $g'(0)=1\neq 0$. Also for any $j$ one 
has
$h_j(0)=0$. By 
the fundamental theorem o Calculus  
$h_1'(x)=g(x)$ and so $h_1$ is at least an $1$--accelerator. Since 
$h^{(j)}_j(x)=g'(x)$ and $h^{(i)}_j(0)= h^{(i-1)}_j(0)=0$ for all $i<j$, 
the result follows.
 \end{proof}
In  
  example 1 of section 
\ref{sec3}
this result is applied to $g(x)=\sin(x)\in\,NEU_\infty$ where
an accelerator of arbitrary order is obtained for $g$.

Another simple way to improve the order of acceleration of an iteration 
function is by composition with itself:
\begin{proposition}\label{prop2}
If $h$ is an accelerator of $\neu$ then the iteration function 
 \begin{equation}\label{eq4}
h_k =\underbrace{h\circ h\circ\ldots\circ h}_{k\,\,\mbox{times}}
\end{equation}
is at least a $k$--accelerator of $\neu$ for $k<m$.
\end{proposition}
\begin{proof}
As $h\in\neu$ is an accelerator then $h(\xast)=\xast$ and $h'(\xast)=0$. By 
the  chain rule  
  $h_2'(x)=h'(h(x))\,h'(x)$ and 
$h_2^{(2)}(x)=h^{(2)}(h(x))\,h'(x)^2+h'(h(x))\,h^{(2)}(x)$. Hence $h'_2 
(\xast)=h_2^{(2)}(\xast)=0$, that is $h_2$ is at least a $2$--accelerator and the result 
follows by mathe\-ma\-tical induction in $k$.  
\end{proof}

 In \cite{graca1,graca3} we introduced and studied other type  of  
 accelerators namely those resulting of a suitable combination of 
 $u\in\, NEU$ with an iteration 
function $v\in \, HYP$. We called these iteration functions {\em combined}.  
\begin{definition}\label{def3}    
Let $u $ be a function with a fixed point $\xast$ and $v\in\, HYP$.
 The following iteration function $h=C(u,v)$,
\begin{equation}\label{eq1}
h(x)=C(u,v)(x)=\frac{v(x)-u(x)\,v'(x)}{1-v'(x)}
\end{equation}
is called a combined iteration function.
\end{definition} 

Note that by  the given definition of combined iteration function it results   that  
 $u$ and $v$ possess the same (generally unknown) fixed point $\xast$. 
 Furthermore, as a  consequence 
of the definition of $h=C(u,v)$ we have 
   $h'(\xast)=0$ 
 when $u\in\, NEU$. So the next proposition.
\begin{proposition}\label{prop3} 
Let $u\in\, NEU_m$ and $v\in\, HYP_j$ with $m\geq 1$ and $j\geq 2$. Then,
  the combined function $h=C(u,v)$ is an 
accelerator of $u$.
\end{proposition}
 
Let us remark  that 
for iteration functions verifying $|v'(\xast)|\neq 1$, that is for 
$v\in\, HYP$, 
a simple accelerator is  
$h(x)=C(x,v)$.   Thus $h$ 
accelerates $HYP$, either the sequence $x_{n+1}=v(x_n)$, for $x_0$ 
sufficiently close to $\xast$,
  is convergent to $\xast$ or not.

Neutral iteration functions satisfying the following hypotheses $(H_1)$ often 
arise in the applications:
 
\begin{equation}\label{nova1}
(H_1)\qquad
\left\{
\begin{array}{ll}
\tilde u(\xast)&=\xast,\\
\tilde u'(\xast)&=1,\\
\tilde u^{(2)}(\xast)&=\alpha\neq 0,\qquad \alpha\in \Rb.
\end{array}
\right.
\end{equation}
Define $  NEU _{2}$ as the set of iteration functions in $NEU_{m}$ for 
which $(H_1)$ holds. For $\bar u\in  NEU_{2} $ consider the iteration 
function
\begin{equation}\label{nova2}
\phi(x)=\bar u(x)-\bar u'(x)+1.
\end{equation}
\begin{proposition}\label{novaA}
The set $ NEU_{2} $ is accelerable. An accelerator is the iteration 
function $h=C(x,\phi)$:
\begin{equation}\label{nova3}
h(x)=\frac{\phi(x)-x \phi'(x)}{1-\phi'(x)},
\end{equation}
where $\phi$ is defined as in {\rm(\ref{nova2})}.
\end{proposition}
\begin{proof}
Let $\bar u\in  NEU_{2} $. Since 
$\phi(\xast)=\xast$ and $\phi'(\xast)=1-\alpha\neq 1$, hence $\phi\in 
HYP_j$  ($\,j\geq 2$) and 
proposition \ref{nova3} holds.
\end{proof}

\noindent
We remark   that when $u\in\neu$ the combined function $C(x,u)$ is not 
 defined at $\xast$. However under the   hypotheses 
 \begin{equation}\label{nova4}
(H)\qquad\left\{ 
\begin{array}{ll}
u(\xast)=&\xast,\\
u'(\xast)=&1,\\
u^{(j)}(\xast)=&0,\qquad2\leq j\leq m-1,\\   
 u^{(m)}(\xast)=&\alpha\neq 0, 
\end{array}
\right.
\end{equation} (which are stronger than $(H_1)$)  
 we can continuously extend $C(x,u)$ in order to get 
$v\in\, HYP$. Namely, for $u\in\, NEU_{m}$, 
\begin{equation}\label{eq6}
 v(x)=\left\{
\begin{array}{ll}
C(x,u)(x)=\displaystyle{\frac{u(x)-x u'(x)}{1-u'(x)}} &x\neq\xast\\
\xast & x=\xast .
\end{array}
\right.  
\end{equation}
is in $HYP$ since $0< v'(\xast) =1-1/m<1$ (see  for instance 
\cite[pp. 99]{Keller1}). So, by proposition \ref{prop3} we obtain the 
following result.
\begin{proposition}\label{prop5}  
Let $u\in \neu$   and   $v$  as  in {\rm(\ref{eq6})}. If $u$ satisfies 
$(H)$ then the iteration function 
 \begin{equation}\label{7}
w(x)=C(x,v)(x)=\displaystyle{\frac{v(x)-x v'(x)}{1-v'(x)}}
\end{equation}
 is an accelerator of $u$.  
\end{proposition}

The map $w$ given in proposition \ref{prop5} is a powerful 
accelerator. Hereafter whenever we use the symbol $w$ we are referring 
the map
  $w$ given in proposition \ref{prop5} and we will call it
the  {\em standard accelerator}  of $u$.
 The iteration function $v$ given by (\ref{eq6}) is   Newton's 
 iteration function $N(x)=x-\psi(x)/\psi'(x)$ for $\psi(x)=x-u(x)$ and   
    $w$ 
is  the  Newton's ite\-ra\-tion function     for $\psi(x)=(x-u)/(1-u'(x))$ (see for instance \cite[pp. 
127]{Traub1}). Hence the final map $w$ is directly obtained as a double 
newtonisation. 
 Of course any other 
iteration function $v\in\,  HYP$ can be combined with $u\in\, NEU$ in order to 
obtain an accelerator.
  The iteration function $v$ given in (\ref{eq6}) is a good choice  but  
  many other possibilities are available. For instance the well known Steffensen's iteration function 
  $v(x)=x-(u(x)-x)^2/(x-2 u(x)+u(u(x)))$ (\cite[pp. 315]{Brezinski1}), 
which is a reformulation of the Aitken's $\Delta^2$ process 
\cite{Aitken1}, 
can be used  or any other ite\-ra\-tion function known to transform a neutral 
fixed point $\xast$ of $u$ into an attracting one.

Between the possible accelerators for a given map the best one one can get 
is an accelerator that transforms any initial point into the fixed point. 
This corresponds to the notion of the kernel of an accelerator.

 \begin{definition}\label{def2}
For $u,  \, v  $ and $w $ as   in proposition \ref{prop5}, 
 the set
$$
ker_w=\left\{
u\in\neu: w(x)=\xast 
\right\}
$$
is called the kernel of $w$. 
\end{definition}
 In examples 3 and 4 of section \ref{sec3}  we show that several iteration 
 functions ge\-ne\-rating 
logarithmic fixed point sequences belong
to the kernel   of the standard accelerator $w$.

A characterization of the  kernel of the standard accelerator in terms of the initial 
iteration  map  $u$ is given by proposition \ref{prop8}. Let us first 
prove a lemma  characterizing  the kernel of $w$ in terms of the 
auxiliary function $v$.   
 \begin{lemma}\label{prop7}
Let $v$ as in {\rm(}\ref{eq6}{\rm)} and $w$ as in {\rm(}\ref{7}{\rm)}. An iteration function $u\in\, NEU$ belongs to the kernel 
of $w$ if and only if
$v$ is the linear function
\begin{equation}\label{12a}
v(x) = (1-v'(\xast))\,\xast + v'(\xast)\, x,
\end{equation}
for $x$ belonging to a neighbourhood of the {\rm(}isolated{\rm)} fixed point $\xast$.
\end{lemma}
\begin{proof}
By definition $u\in Ker_w$ if and only if $w$ is the constant function 
$w(x)=\xast$, for all $x$. So, $w'(x)=0\,\,\forall x$. Differentiating 
$w(x)$ given by (\ref{7}) 
we have,
$$
\begin{array}{ll}
w'(x)=0&\Longleftrightarrow\,\, v''(x)\left( v(x)-x
\right)=0\\
& \Longleftrightarrow v''(x)=0\,\,\mbox{or}\,\, v(x)=x.
\end{array}
$$
As by hypothesis $v'(\xast)\neq 1$ then the case $v(x)=x$ is excluded. So, $ 
w'(x)=0$ iff $v''(x)=0$. That is, $v(x)$ is of the form $v(x)=a\,x+b$. As 
$v(\xast)=\xast$ and $v'(\xast) $ exists then $a=v'(\xast)$ and 
$b=(1-v'(\xast)) \xast$.
\end{proof}

\begin{proposition}\label{prop8}
An iteration function $u\in NEU_m$ satisfying $(H)$ belongs to the 
kernel  of the standard map $w$  if and only if $u$ has the form
\begin{equation}\label{A}
u(x)=x+\alpha\,  (\xast -x)^\beta ,
\end{equation}
for some constants $\alpha$ and $\beta$   with $\beta \neq 0$.
\end{proposition}
\begin{proof}
Let $u\in NEU_m$. By lemma \ref{prop7} there are constants $a\neq 1$ and 
$b=(1-a)\,\xast$ such that 
\begin{equation}\label{B}
v(x)=a\,x+b
\end{equation}
where
\begin{equation}\label{C}
v(x)=\frac{u(x)-x\,u'(x)}{1-u'(x)}\in HYP.
\end{equation}
Equating (\ref{B}) and (\ref{C})
 we have that $u$ belongs to the kernel of $w$ if and only if it is a 
 solution  of the
following first order   differential equation
\begin{equation}\label{D}
(b+(a-1)\,x)\,u'(x)+u(x)=ax+b.
\end{equation}
The general solution  of (\ref{D}) is
\begin{equation}\label{E}
 u(x)=x+c\,(b-(1-a)x)^{\frac{1}{1-a}},\,\, \forall c\in \Rb.
\end{equation}
The solution $u$ given by (\ref{E}) satisfies  $u(\xast)=\xast$ and 
$u'(\xast)=1$ for an arbitrary constant $c$ since $b=(1-a)\,\xast$. So, 
taking $\alpha= c\,(1-a)^{\frac{1}{1-a}}$ and $\beta=1/(1-a)$ the result follows.
\end{proof}
 
\section{The standard accelerator and logarithmic fixed point sequences}\label{sub1}

As remarked by Delahaye in \cite[pp. 181]{del}, "logari\-thmi\-cally convergent sequences 
are difficult to acce\-le\-rate: it is not possible to accelerate all of them 
with only one transformation" (that is $Log$ is not accelerable). However the 
scenario is quite different for fixed point sequences.

Let us define $hyp_j $ as being the set of sequences generated by  
iteration functions in $\hyp$ and   $Log_m \subset Log$  defined as 
follows:
 \begin{equation}\label{10}
Log_m =\left\{
(x_n):  
x_{n+1}=u(x_n),  \,\,
x_0\in {\cal D},\,u\in\neu\,\,
 \mbox{and}\,\, u\,\, \mbox{satisfies}\,\, 
 (H)  
\right\}. 
\end{equation}

We show in the next proposition that the set 
  $FIX=hyp_j \cup \,Log_m $,   
  of  fixed point sequences 
 can be acce\-lerated using a single   transformation 
namely the one associated to our standard map $w$.  

We say that a sequence to sequence transformation $W:$ $\,\xn\rightarrow$ $ (w_n)$
accelerates the convergence of $\xn$  
if
  $(w_n)$ $\setan \xast$ and $\limn (w_n-\xast)/(x_{n}-\xast)=0$.
 
Let $\xn\in FIX$ and $W:\,\xn\rightarrow (w_n)$ be the sequence to sequence 
transformation 
such that 
  $w_n=(w(x_{n-1})), \,\,n=1,2,\dots$  where $v$ and $w$ are   iteration functions  
  as in 
proposition  \ref{prop5}.
\begin{proposition}\label{prop6}
 For any natural numbers $m$ and $j$, and for $x_0$ chosen sufficiently close to 
 $\xast$, the   set   $FIX$
 is accelerable by the sequence to sequence transformation $W$.
 \end{proposition}
 \begin{proof}
If $\xn\in FIX$ then either  $\xn\in hyp_j$  or $\xn\in 
Log_m$ (note that these two sets do  not intersect).
 By definition of $FIX$ the sequence $\xn$ is generated by 
an iteration function $u$ i.e. $x_n=u(x_{n-1})$. If $u\in\hyp$ then both 
$v=C(x,u)$ and $w=C(x,v)$ verify $v'(\xast)=0=w'(\xast)$ since they are 
combined iteration functions.

If $\xn\in Log_m$  and
  $(H)$ holds then
  proposition \ref{prop5} yields also $w'(\xast)=0$. So, either for $\xn$ 
  in $Log_m$ or $hyp_j$,    the chain 
rule  and the mean value theo\-rem applied to the function 
$F=w\circ v\circ u$ enable us to claim that there exists a neighbourhood 
$\Omega$ of $\xast$ such that for any $x_0\in\Omega$ the sequence $(w_n)$ 
converges to $\xast$ and   $\limn (w_n-\xast)/(x_{n}-\xast)=0$. That is
  $W$ accelerates $FIX$.
\end{proof}

\section{Examples}\label{sec3}

 \noindent
  {\it Example} 1. The   iteration function $g(x)=\sin x$ is frequently used as a 
test function for 
assessing the quali\-ty of sequence to sequence transformations 
 \cite[pp. 325]{Brezinski1}. Popular 
sequence transformations such as the $\epsilon$   and $\rho$ 
   algorithms (see respectively  \cite{wynn1} and \cite{wynn2}) , 
Richardson extra\-polation  \cite[pp. 181]{gau} , iterated Aitken $\Delta^2$  
 (\cite{Aitken1},\cite[pp. 229]{del} and for instance\cite{wen}), Overholt 
process  \cite{over}  and Levin 
transformation  \cite{lev}  are unable to accelerate $g$. The $\Theta$ 
algorithm \cite{Brezinski2}
  produces satisfactory numerical results though it is only a 
$1$--accelerator for the sequence $x_{n+1}=g(x_n)$.

We show below that there are $j$--accelerators for $g$   of 
the type given in \eqref{eq3}. Our standard map $w$ is also compared with 
Aitken's $\Delta^{2}$ process and the $\Theta_2$ procedure. As can be seen 
in table~\ref{tabelaA1} the map $w$ performs much better than the 
referred processes.

\begin{itemize}
\item[a)]
 The map $g(x)=\sin x$ has  $j$--accelerators,   with $j>1$, of the type  \eqref{eq3}. Indeed  
$$
g(x)=\sin x\,\,\in  NEU_\infty,\quad g(0)=0 
$$
 so  by proposition  \ref{prop1} the function $g$ is $\infty$--accelerable.
The first  $h_j$--accelerators of $g$ are
$$
\begin{array}{l}
h_1(x)=\int_0^x \sin t\, dt=1-\cos x\,\, \\
{}\\
h_2(x)=\int_0^x h_1(t)\, dt=x-\sin x\,\,  \\
{}\\
h_3(x)=\int_0^x h_1(t)\, dt=\frac{x^2}{2} -1+\cos x . \\
\end{array}
$$
 
 It is easily deduced an 
explicit formula for $h_j$, ($\,\,j\geq 1$) :  
$$
h_j(x)=\left\{
\begin{array}{l}
(-1)^ { \frac{j-1}{2}}  \, 
q_{j-1}(x)+(-1)^ { \frac{j+1}{2}}  \cos 
x, 
\qquad  \mbox{for}\,\,j\,\mbox{odd}\\
\\
(-1)^{\frac{j-2}{2}} \, 
p_{j-1}(x)+(-1)^{ \frac{j}{2}} \sin x,\quad 
\qquad \mbox{for}\,\,j\,\mbox{even}.\\
\end{array}
\right. 
$$
where $q_j(x)$ and $p_j(x)$ are the Maclaurin polynomial of 
degree $j$,  respectively of $\cos x$ and $\sin x$. The iteration function $h_j$ is a 
$j$--accelerator of $g$.
\item[b)]
For $g(x)=\sin(x)$ the maps $v=C(x,g)$ and $w=C(x,v)$ are respectively,
$$\begin{array}{l}
\displaystyle{v(x)=\frac{x\,\cos(x)-\sin(x)}{\cos(x)-1}},\\
\\
\displaystyle{w(x)=\frac{(x^2-1)\,\sin(x)+\cos(x)\,(x+\sin(x))-x}{2\,\cos(x)+x\,\sin(x)-2}}.
\end{array}
$$
In table \ref{tabelaA1} 
is displayed the first $4$ iterations for the maps $g$, $v$ and $w$ with 
starting point $x_0=3.0$. These results are compared with the Aitken's $\Delta^2$ 
process and the $\Theta_2$ procedure.

\begin{table}[ht] 
\caption{  \label{tabelaA1} Iterations for the maps $g$, $v$, $w$ with starting 
point $x_0 =3.0$. Comparison with $\Delta^2 $ 
and $\Theta_2 $ processes. }  
$$
\begin{array}{|l|l|l|l|l|}
\hline
\qquad g&\qquad v&\qquad w&\qquad \Delta^2&\qquad \Theta_2\\
\hline
3.0&3.0&3.0& &\\
0.14112&1.56337&1.40041&0.140652& \\
0.140652&0.995758&0.173163&0.0938926&0.141125\\
0.140189&0.652467&0.000345858&0.0935825&- 0.000754788\\
0.13973&0.431844&7.30548\times 10^{-13}& &\\
\hline
\end{array}
$$
\end{table}

\end{itemize}


\bigbreak
 \noindent
 {\it Example} 2.
The well known logistic map \cite[pp. 2]{stuart1}
$$
\bar u(x)=a\,x\,(1-x) 
$$
has the neutral fixed point $\xast=0$ when $a=1$. Since $\bar u'(0)=1$ and $\bar 
u^{(2)}(0)=-2\neq 0$ then $\bar u\in  NEU_{2}$. So,  by proposition 
\ref{novaA}, for
$$
\phi(x)=\bar u(x)-\bar u'(x)+1=3-x^2,
$$
the iteration map
$$
h(x)=C(x,\phi)(x)=\frac{x^2}{2\,(x-1)},\quad x\neq 1,
$$
is an accelerator for $\bar u$. The table \ref{tabelanova} displays some 
iterations respectively for  $y_n=\bar u (y_{n-1}) $ and   
$x_n=h(x_{n-1}) $ taking $y_0=x_0=0.5$ as starting point.

When $a\neq 0,1$, the fixed point is 
$\xast=\displaystyle{\frac{a-1}{a}}$ and $\bar 
u\in HYP$ since $\bar u'(\xast)=2-a\neq 1$. So, by proposition 
\ref{prop3}  the iteration map $h=C(x,\bar u)$ is an accelerator for 
$\bar u$:
$$
h(x)=\frac{\bar u(x)-x\,\bar u'(x)}{1-\bar u'(x)}=\frac{a x^2}{1-a+2a 
x},\qquad x\neq \frac{a-1}{2a}.
$$

\begin{table}[ht] 
\caption{  \label{tabelanova} $y_n=\bar u (y_{n-1}) $, $x_n=h(x_{n-1}) $ }  
$$
\begin{array}{|l|l|}
\hline
\qquad y_n&\qquad x_n\\
\hline
0.5&0.5\\
0.25&-0.25\\
0.1875&-0.025\\
0.1523\ldots&-0.000304\ldots \\
\hline
\end{array}
$$
\end{table}


\bigbreak
\noindent
 {\it Example 3.}  
The following  iteration function $g$ (cf. \cite{fdil}) belongs to the kernel of 
 the standard map $w$: 
 $$
g(x)=x+(x-1)^{3/2},\quad \xast=1 . 
$$
Note that $g\in NEU_1$ and $g'(1)=1$. We have,
$$
v(x) =\displaystyle{\frac{g(x)-xg'(x)}{1-g'(x)}} =
 \displaystyle{\frac{x+2}{3}}=C(x,g), 
$$
and $ w(x)=C(x,v)(x)=1=\xast$.  In fact as
 the iteration function $v$ has the form
$$
  v(x)= (1- v'(\xast))\xast+v'(\xast)\,x,
$$
then, by  lemma \ref{prop7}, $g$ belongs to the kernel of $w$.

 Proposition \ref{prop8} can be directly applied to the computation of the 
 multiple zero of a polynomial of the type
$$
p(x)=\alpha\, (x-\zeta)^r,
$$ 
where $\alpha\neq 0$ and $r\geq 1$. Indeed
 $p$ has a zero $\zeta$ of multiplicity $r$ and
$$
u(x)=x-p(x)=x-\alpha\,(x-\zeta)^r
$$
has the form (\ref{A}) of proposition \ref{prop8}, that is $u$ belongs to $Ker_w$.

More generally, for $\alpha\neq 0 $ and $\, r\in\Rb,\,\, r> 1$ 
consider the set
 $$
{\cal P}=\left\{
g \in NEU_\infty :\,
 g(x)=x+\alpha (x-\xast)^r 
\right\}.
$$

Computing $v=C(x,g)$, for $g\in\, {\cal P}$, one has
$$
v(x)=(1-\frac{1}{r}) x+\frac{1}{r} \xast.
$$
Thus $v$ has the form (\ref{12a}) of lemma \ref{prop7} and so     any iteration 
function in
  $ {\cal P}$ belongs to $Ker_w$.


\bigbreak
\noindent
 {\it Example 4.}   In \cite{sab1,sab2,sedogbo} several sequence to sequence transformations have been 
developed and 
tested in order to find an accelerator for logarithmic sequences 
generated by iteration functions belonging to subsets of the following set 
$$
{\cal S}=\left\{
\begin{array}{l}
g \in NEU_\infty  :
\, g(x)=x+\sum_{i=1}^{\infty}\alpha_i (x-\xast)^{r+i} \\
 \end{array}
\right\},
$$
where $\alpha_1\neq 0$, $r\in\Rb,\,\,r\geq 1$.  
 
  The  referred
authors   hardly found an accelerator for some subsets of ${\cal S}$ in spite 
of considering only  
  positive integers   $r$ and $\alpha_1<0$.  

The set ${\cal S}$ is accelerable by our standard 
map $w$. Indeed, since
$$
g'(x)=1+\alpha_1 (r+1) (x-\xast)^r+q_r(x),
$$
where $q_r(x)$ is a power  series   of $(x-\xast)$ whose first 
term has degree $ r+1$  we get 
$g\in NEU_\infty  $, $g^{(j)}(\xast)=0,\quad 2\leq j\leq r$ and  
 $\,g^{(r)}(\xast)\neq 0$. So  proposition \ref{prop5} holds.


\bigbreak
\noindent
 {\it Example 5.}

Although the combined iteration functions have been defined for  real 
functions they can 
be generalized to  complex maps. The extension of the results obtained 
in this paper to higher dimensional maps   will be object of another work.

In order to illustrate the behaviour of the standard (complex) accelerator 
we chose the analytic function
\begin{equation}\label{6.1}
f(z)=z^2\,(z-2)^2\,(e^{2z}\,\cos z+z^3-1-\sin z)
\end{equation}
studied by Kravanja and Van Barel \cite[Example 4]{kravanja1}.
These authors have used this function to test their  algorithms 
for the simultaneous computation of zeros of analytic functions. As stressed 
in   
the referred work the accuracy
of the numerical results obtained by such algorithms diminishes as the 
multiplicity of a zero increases.  The best numerical results obtained in 
\cite{kravanja1} for the multiple roots $\zeta_1=0$   and
$\zeta_2=2$ (both having only multiplicity $2$) do not exceed $7$ correct 
digits.

  Our combined maps do  not need the 
knowledge of the multiplicity of a zero and produce accurate numerical 
results in a few iterations as shown in   table \ref{tabela1} where, for 
instance, the 
 multiple zero $\zeta_2$ is accurately computed in 5 iterations of the 
 standard map. However, the price to pay with our methods is the need for
an initial approximation for each zero (chosen, say, to be within an euclidean distance 
about $0.1$ from a zero of $f$).

The basic iteration function to be used is
$
u(z)=z-f(z)
$,
which can lead to complex divergent sequences $y_n=u(y_{n-1})$, with 
$y_0\in \Cb$ 
taken as initial approximation to a zero of $f$ (see table \ref{tabela1}). The sequence obtained by 
iteration of the standard complex map is denoted by  $(z_n)$ with 
$z_0=y_0$.  The numerical example presented here has been 
 programmed in {\sl Mathematica} \cite{wolfram1} and  the computations 
were carried out with $16$ digits of precision.
Due the lack of space in table \ref{tabela1} the imaginary part of each iterate has been 
partially truncated.
\begin{table}[ht] 
\caption{  \label{tabela1} First $6$ iterations for $u(z)$ and $w(z)$, with 
$z_0=1.9+0.1\, i$.}  
$$
\begin{array}{|l|l|}
\hline
\qquad\qquad\qquad y_n=u(y_{n-1})&\qquad\qquad\qquad  z_n=w(z_{n-1})\\
\hline
 1.900000000000000+0.10000000\, i&1.900000000000000 + 
 0.1000000\, i\\
 2.391422135261736 - 0.4699667\, i&2.033556020548597 + 
 0.0804529 \, i\\
 -190.6272479365824 + 83.78040\, i&2.010056510555553 - 0.01717596
 \,i\\
 -1.078985533\times 10^{45}+ 2.057\times 10^{45}\, i& 2.000502552323976 + 0.0010266 
\, i\\
      {"Indeterminate"}&2.000002378186929  - 3.083\times 10^{-6}\, i\\
      {"Indeterminate"}& 1.999999999946048+0.\times 10^{-16}\, i \\
\hline
\end{array}
$$
\end{table}



\newpage 
\begin{thebibliography}{99}


\bibitem{Aitken1} A. C. Aitken,  On {B}ernoulli's numerical solution of algebraic
equations, {\it Proc. Roy. Soc. Edinb.}, 46, 1926,  pp.~289--305.

 

\vspace{-0.4pc}
\bibitem{Brezinski1} C. Brezinski and M. Redivo Zaglia,
{\it Extrapolation Methods Theory and Practice},
North--Holland, Amsterdam, 1991.

\vspace{-0.4pc}
\bibitem{Brezinski2} C. Brezinski,  Acc\'{e}l\'{e}ration de suites \`{a} 
convergence logarithmique, {\it C. R. Acad. Sci. Paris}, 273 A, 1971,  pp.~727--730.

\vspace{-0.4pc}
  \bibitem{fdil} A. Fdil, A new method for solving a nonli\-near equation 
  with error estimations,
 {\it Appl. Numer. Math.}, 21,
1996, pp.~417--429.

\vspace{-0.4pc}
\bibitem{del} J.~P. Delahaye,{\it Sequence Transformations},
 Springer--Verlag, 1988.

\vspace{-0.4pc}
\bibitem{delgb} J.~P. Delahaye and B. Germain--Bonne, The set of 
logarithmically convergent sequences cannot be accelerated,
 {\it SIAM J. Numer. Anal.}, 19, 1982,  pp.~840--844.
   
 

\vspace{-0.4pc}
 \bibitem{gau} W. Gautschi, {\it Numerical Analysis an Introduction}, 
  Birkh\"{a}user, Boston, 1997.

 \vspace{-0.4pc}
 \bibitem{graca1} M.~M. Gra\c{c}a, Nonhyperbolic fixed 
 points and flat iteration functions, {\it (to appear)  Experimental 
 Mathematics} 11:1, 2002.

  
 
\bibitem{graca3} M.~M. Gra\c{c}a, Simple accelerators for logarithmic 
fixed point sequences, {\it Problems in Applied Mathematics and 
Computational Intelligence, N. Mastorakis Editor}, World Scientific and 
Engineering Society Press, 2001, pp. 82--86.
 

\vspace{-0.4pc}
\bibitem{Hol} R.~A. Holmgren,  {\it A First Course in Discrete Dynamical 
Systems}, 
  Springer, N.Y., 1996. 

\vspace{-0.4pc} 
\bibitem{Keller1}  E. Isaacson \& H. B. Keller,
{\it Analysis of Numerical Methods},
 John Wiley \& Sons,
 N. Y., 1966.

 \vspace{-0.4pc}
\bibitem{kravanja1} P.~Kravanja and M.~Van Barel, A derivative--free 
algorithm for computing zeros of analytic functions, {\it Computing}, 63, 
1999, pp.~69--91.
 
 
\vspace{-0.4pc}
\bibitem{lev} D. Levin,  Development of non--linear transformations for 
improving convergence of sequences, {\it Int. J. Comput. Math.}, B3, 1973, 
 pp.~371--388.
 
\vspace{-0.4pc}
\bibitem{over} K.~J. Overholt,  Extended Aitken acceleration, 
  {\it BIT}, 5, 1965, pp.~122--132.

\vspace{-0.4pc}
\bibitem{sab1} P. Sablonni\`{e}re,   Convergence acceleration of 
logarithmic fixed point sequences , 
  {\it J. Comput. Appl. Math.}, 19,  1987 , pp.~55--60.

\vspace{-0.4pc} 
\bibitem{sab2} P. Sablonni\`{e}re,   Comparison of four algorithms 
accelerating the convergence of some logarithmic fixed point sequences , 
  {\it Numer. Algorithms 1}, 1991,  pp.~177--198.

 \vspace{-0.4pc}
\bibitem{sedogbo} G.~A. Sedogbo,   Convergence acceleration of some 
logarithmic sequences, 
  {\it J. Comp. Appl. Math.}, 1990,  pp.~253--260.

\vspace{-0.4pc}
\bibitem{stuart1} A.~M. Stuart and A.~R. Humphries,  {\it  Dynamical 
Systems and Numerical Analysis}, 
  Cambridge University Press, Cambridge, 1998. 


\vspace{-0.4pc} 
\bibitem{Traub1} J.~F. Traub, {\it Iterative Methods for the Solution of Equations}, 
  Prentice-Hall, N.J., 1964.

\vspace{-0.4pc} 
\bibitem{wen} E.~J. Weniger,  On the derivation of ite\-ra\-ted sequence transformations for 
the acce\-le\-ration of convergence and the summation of divergent series,
 {\it Comput. Phys. Commun.}, 64, 
 1991, pp.~19--45.

\vspace{-0.4pc}
\bibitem{wynn1} P. Wynn,  On a device for computing the $e_m(S_n)$ 
transformation, 
 {\it MTAC}, 10, 
 1956, pp.~91--96.

\vspace{-0.4pc}
\bibitem{wynn2} P. Wynn,  On a procrustean technique for the numerical 
transformation of slowly convergent sequences and series,
 {\it Proc. Cambridge Phil. Soc.}, 52, 
 1956, pp.~663--671.
 
 \bibitem{wolfram1} S. Wolfram, {\it The Mathematica Book}, 
   Third ed., Cambridge University Press, 1996.

 \end{thebibliography}
\end{document}